\magnification= 1150
\hoffset= .22in
\voffset=-.2in
\input amstex
\documentstyle{amsppt}
\input xypic
\NoBlackBoxes
\NoRunningHeads



\define\incl{\subseteq}

\define\wt{\widetilde}                                %
\define\ov{\overline}                                 %
\define\tagit#1{\tag "(#1)"}
\define\spec#1{\underline{\underline{#1}}}

\define\Wit2#1{W_*(#1) \tsize{[\frac 12]}}
\define\Wt2#1{\underline{\underline{W}}(#1)\tsize{[\frac12]}}


\define\prf{\demo{\underbar{Proof}}}
\define\endpf{\enddemo}
\define\dfn#1{\definition{\bf\underbar{Definition #1}}}
\define\enddef{\enddefinition}


\define\lamg#1{\sum \lambda_{#1i}g_{#1i}}
\define\insqx#1{\{x_j\, |\, x_j \in \BbbC[#1]\}}


\define\cgpph#1{C^*(#1;\varphi)}
\define\cgpal#1{C^*(#1;\{\varphi_\alpha\})}            %
\define\Hi#1{H^\infty(#1)}                          %

\define\hil#1{H^{\infty}_L(#1)}                               %
\define\hwtil#1{H^{\infty}_L(#1)}                               %
\define\hpil#1{H^{\infty}_{L'}(#1)}
\define\gpal#1{\BbbC[#1]}

\define\yinn1{I^n_{n+1}[\Gamma]}                               %

\define\xjjn1{J^j_{n+1}(\Gamma)}
\define\jj1n1{J^{j+1}_{n+1}(\Gamma)}

\define\yistnn1{I^{n,*}_{n+1}(\Gamma)}

\define\xicstj1n1{I^{j+1,*}_{n+1}(\Gamma)}
\define\xicstjn1{I^{j,*}_{n+1}(\Gamma)}

\define\istzer1{I^{0,*}_1(\Gamma)}

\define\ixnn1A{I^n_{n+1}(A)}


                   %
                          %

                        %
                        %
                                      %
                                 %
                                     %
                                    %
                                    %

                               %
                                     %
                           %


\define\BbbC{\Bbb C}                                  %
\define\BbbQ{\Bbb Q}                                  %
\define\BbbX{\Bbb X}                                  %
\define\BbbZ{\Bbb Z}                                  %
\define\X00{\BbbX^0_0(\Gamma)}
\define\xS00{S^0_0(\Gamma)}
\define\bS0{\underline S_0(\Gamma)}


\define\gm{\Gamma}                                    %

\define\g0{\Gamma_0}                                  %
\define\rhbpi#1{H_{*+#1}(B\pi\,;\,\BbbQ)}
\define\gl#1{GL(#1)}
\define\xgl1#1{GL_1(#1)}
\define\ygln#1{GL_n(#1)}
\define\mnn#1{M_n(#1)}


                                       %
                                       %

                             %
                           %
                           %
                           %
                           %
\define\x0i{x_{0i}}
\define\xtil0i{\widetilde x_{0i}}
\define\zee#1{\BbbZ[#1]}


\define\sumoverg1{\sum\limits_{i\,|\,\overline g_i=\overline g}}
\define\intjn{\bigcap^j_{k=0} \ker\,(\partial_k :}
\define\xintjn2{\bigcap^j_{k=0} \ker\,(\partial_k \times\partial_k :}

\define\prn1{\prod^1_{i = n}}
\define\pr0n{\prod^n_{i = 0}}
\define\prj0j{\prod^j_{i = 0}}


\define\surj{\twoheadrightarrow}                      %
\define\inj{\rightarrowtail}                          %

\define\twobytwo#1#2#3#4{
\left[\vcenter{\offinterlineskip                            %
\halign{\strut##&\hfil##\hfil&\ \vrule##\ &\hfil##\hfil\cr  %
&#1&&#2\cr                                                  %
\noalign{\hrule}                                            %
&#3&&#4\cr}}\right]}                                        %
                                          %
                                       %
                       %
                       %
%
                                        %
                                 %
                                      %

 %
                         %

\document
\baselineskip 20pt


\vskip.5in
\topmatter
\title
Filtrations of simplicial functors and the Novikov Conjecture
\endtitle
\vskip.2in
\author
C. Ogle
\endauthor
\affil
The Ohio State University
\endaffil
\date Oct., 2004; revised Feb. 2005
\enddate
\abstract We show that the Strong Novikov Conjecture for the maximal $C^*$-algebra
$C^*(\pi)$ of a discrete group $\pi$ is equivalent to a statement in topological
$K$-theory for which the corresponding statement in algebraic $K$-theory is always true. We also show that for any group $\pi$, rational injectivity of the full assembly map for $K_*^t(C^*(\pi))$ follows from rational injectivity of the restricted assembly map. 
\endabstract
\toc
\widestnumber\subhead{3.2.1}
\subhead 0. Introduction
\endsubhead
\subhead 1. Simplicial $C^*$-algebras
\endsubhead
\subhead 2. Filtering homotopy groups
\endsubhead
\subhead 3. A basic example
\endsubhead
\subhead 4. Filtering the assembly map
\endsubhead
\subhead 5. Simplicial $C^{i}$-algebras and the Strong Novikov Conjecture 
\endsubhead
\subhead 6. Some additional remarks
\endsubhead
\endtoc
\address
Dept. of Mathematics,
The Ohio State University
\endaddress
\email
ogle\@math.ohio-state.edu
\endemail
\endtopmatter


\newpage
\head 0 Introduction
\endhead
\vskip.2in

Let $\pi$ be a countable discrete group. The Novikov conjecture for $\pi$ is
equivalent to the assertion that the assembly map for the Witt groups of $\Bbb Z[\pi]$,
$$
H_*(B\pi;\underline{\underline{W}}(\Bbb Z))\to W_*(\Bbb Z[\pi])
\tag0.1
$$
is rationally injective (here
$H_*(B\pi;\underline{\underline{W}}(\Bbb Z))$ denotes the homology
of $B\pi$ with coefficients in the spectrum
$\underline{\underline{W}}(\Bbb Z)$). The Novikov Conjecture
states the Novikov Conjecture for $\pi$ is true for all discrete
groups $\pi$ [N]. Let $C^*(\pi)$ denote the maximal $C^*$ algebra
of $\pi$. As with Witt theory, there is an assembly map
$$
KU_*(B\pi) = H_*(B\pi;\underline{\underline{K^t}}(\Bbb C))\to K^t_*(C^*(\pi))
\tag0.2
$$
where $\underline{\underline{K^t}}$ denotes topological $K$-theory
spectrum, and $K_*^t(_-)$ its homotopy groups. The Strong Novikov
Conjecture (SNC) for $\pi$ asserts that the above assembly map for
$C^*(\pi)$ is rationally injective [Ka]. The Strong Novikov
Conjecture is the statement that SNC is true for all discrete
groups $\pi$. Variants and generalizations of these conjectures
exist for both topological and algebraic $K$-theory ([BC], [BCH],
[C1], [C2], [DL], [FJ], [GV], [HR], [Y1]), and have been verified
for many types of discrete groups ([Ka], [Mi], [KS], [BHM], [CM],
[CGM], [Y2] ).

After tensoring with $\Bbb Q$, the assembly maps in both (0.1) and (0.2)
can be alternatively written as
$$
\gathered
\oplus_{n\in\Bbb Z} H_{*- 4n}(B\pi;\Bbb Q)\to W_*(\Bbb Z[\pi])\otimes\Bbb Q\\
\oplus_{n\in\Bbb Z} H_{*- 2n}(B\pi;\Bbb Q)\to K^t_*(C^*(\pi))\otimes\Bbb Q
\endgathered
\tag0.3
$$
In both cases, there is a \underbar{connective assembly map} 
$$
\gathered
\oplus_{n\in\Bbb N}  H_{*-  4n}(B\pi;\Bbb Q)\to W_*(\Bbb Z[\pi])\otimes\Bbb Q\\
ku_*(B\pi)\otimes\Bbb Q = \oplus_{n\in\Bbb N}  H_{*- 2n}(B\pi;\Bbb Q)\to K^t_*(C^*(\pi))\otimes\Bbb Q
\endgathered
\tag0.4
$$
given by the restriction of the assembly maps in (0.3) to summands indexed by $n\ge 0$, and a \underbar{restricted assembly map} 
$$
\gathered
H_{*}(B\pi;\Bbb Q)\to W_*(\Bbb Z[\pi])\otimes\Bbb Q\\
H_{*}(B\pi;\Bbb Q)\to K^t_*(C^*(\pi))\otimes\Bbb Q
\endgathered
\tag0.5
$$
given by the further restriction to the summand corresponding to $n=0$. From the periodicity isomorphisms $W_*(R)\cong W_{*+4}(R)$ and $K_*^t(A)\cong K_{*+2}^t(A)$ ($R$ a discrete ring with involution, $A$ a complex Banach algebra), one easily sees that the assembly map of either (0.1) or (0.2) is rationally injective if and only if its corresponding connective assembly map is rationally injective. However, it is not obvious that one can restrict any further. Our first result is (cf. Corollary 4.5):

\proclaim{\bf\underbar{Theorem A}} For a given discrete group $\pi$ the assembly map in (0.2) is rationally injective if and only if the corresponding restricted assembly map in (0.5) is injective.
\endproclaim

Although we do not give an explicit proof in this paper, the techniques used in verifying Theorem A hold equally well for $W_*(\Bbb Z[\pi])$, yielding the analogue of Theorem A for the Witt groups of $\Bbb Z[\pi]$.

Our second result consists of i) a reformulation of
SNC in topological $K$-theory, and ii) a proof of the
\underbar{algebraic} $K$-theoretic analogue of the reduction in
i). Given an augmented free simplicial resolution
$\Gamma\hskip-.03in .^+$ of a discrete group $\pi$, we define a
sequence of $\Bbb C$-algebras $\wt
C^i_{n-1}(\Gamma\hskip-.03in .^+),\,\,n\ge 1$, abelian groups
$A_n$, and maps
$$
\gathered
A_n\surj H_n(B\pi;\Bbb Q)\\
 A_n\to K_1^a(\wt
C^i_{n-1}(\Gamma\hskip-.03in .^+))\otimes\Bbb Q\to K_1^t(\wt
C^i_{n-1}(\Gamma\hskip-.03in .^+))\otimes\Bbb Q
\endgathered
$$
where $K_1^a(\wt C^i_{n-1}(\Gamma\hskip-.03in .^+))$ resp.
$K_1^t(\wt C^i_{n-1}(\Gamma\hskip-.03in .^+))$ denotes the the
first algebraic resp. topological $K$-group of $\wt
C^i_{n-1}(\Gamma\hskip-.03in .^+)$, equipped with the fine
topology, and the second map is induced by the natural
transformation from algebraic to topological $K$-theory.
\vskip.2in

\underbar{Statement (A,n)} Each $y\in A_n$ mapping to a non-zero
element under the projection to $H_n(B\pi;\Bbb Q)$ maps to a
non-zero element in $K_1^a(\wt C^i_{n-1}(\Gamma\hskip-.03in
.^+))\otimes\Bbb Q$. \vskip.1in

\underbar{Statement (T,n)} Each $y\in A_n$ mapping to a non-zero
element under the projection to $H_n(B\pi;\Bbb Q)$ maps to a
non-zero element in $K_1^t(\wt C^i_{n-1}(\Gamma\hskip-.03in
.^+))\otimes\Bbb Q$. \vskip.1in

\proclaim{\bf\underbar{Theorem B}} Statement (A,n) is true for all
$n\ge 1$. Moreover, the SNC is true for $\pi$ if and only if
Statement (T,n) is true for all $n\ge 1$.
\endproclaim

The truth of either statement does not depend on choice of free
resolution. Statement (A,n) is true for all free simplicial
resolutions, and if Statement (T,n) is true for one free
simplicial resolution of $\pi$, it is true for all others.

The paper consists of six sections. In section one, we prove some
relevant properties of simplicial $C^*$-algebras, and show how
they can be used to define a filtration on $\pi_0$ (which is a
$C^*$-algebra). In section two we provide an alternative method for
filtering the homotopy (resp. homology) groups of $F(-1)$, where
$F$ is an augmented simplicial object in a suitable category of
spectra, spaces or chain complexes. In section three we identify the
filtration in a basic case, used later in section five. In section four,
we show that two apriori distinct filtrations defined on
$K^t_*(C^*(\pi))$, using sections one and two, in fact agree. In
section 5, we use this to prove Theorem 1. Section six
contains some additional results on the $C^i$-algebras introduced in section five.

This paper is based on an approach to the Novikov
conjecture investigated by the author a number of years ago. We
would like to thank P. Baum for numerous conversations early on
which were helpful in the development of this approach. Also, we
would like to thank D. Burghelea for his continuous support.
\vskip.1in


\newpage

\head 1 Simplicial $C^*$-algebras
\endhead
\vskip.2in

We begin by recalling some basic definitions. The simplicial category $\Delta$ is the
category whose set of objects are the totally ordered sets
$$
\{\underline n = (0<1<2<\dots <n)\} _{n\ge 0}
$$
and whose morphisms are set maps which preserve the ordering. The augmented simplicial
category $\Delta_+$ is formed from $\Delta$ by adjoining an initial object labeled
$\underline{-1}$. A simplicial object in a category $\underline C$ is a contravariant
functor $F:\Delta\to\underline C$, or, equivalently, a covariant functor $\Delta^{op}
\to\underline C$, where $\Delta^{op}$ denotes the opposite category of $\Delta$. Similary,
an augmented simplicial object in $\underline C$ is a covariant functor $F:\Delta^{op}_+
\to\underline C$. A simplicial object is typically represented as a sequence of objects
$\{F(\underline n)\}_{n\ge 0}$ together with face maps $\partial_i: F(\underline n)\to
F(\underline {n-1})$, $0\le i\le n$ and degeneracy maps $s_j: F(\underline n)
\to F(\underline{n+1})$, $0\le i\le n$ satisfying the simplicial identities, while an
augmented simplicial object includes an additional object $F(\underline {-1})$ together
with an \underbar{augmentation map} $\varepsilon: F(\underline 0)\to F(\underline {-1})$
equalizing $\partial_0$ and $\partial_1$.

Suppose $G\hskip-.015in .$ is a simplicial group. For $n \ge 1$ let
$G_n^k = \bigcap_{i=0}^k ker(\partial_i : G_n\to G_{n-1})$. Let $G_0^0 = \partial_1(G^0_1)$.
For notational convenience, we also set $G^{-1}_n = G_n$. Then for each $-1\le k < n$ there
is a split-exact sequence
$$
G^{k+1}_{n+1}\inj G^k_{n+1}\surj G^k_n
\tag1.1
$$
with $G^k_{n+1}\surj G^k_n$ induced by $\partial_{k+1}$ and splitting $G^k_n\inj G^k_{n+1}$
induced by $s_{k+1}$. When $k=n$ there are closed sequences
$$
\gathered
G^{n+1}_{n+1}\inj G^n_{n+1}\overset\partial_{n+1}\to\longrightarrow G^n_n\\
G^0_0\inj G^{-1}_0 = G_0\surj G_0/G^0_0
\endgathered
\tag1.2
$$

The (combinatorial) \underbar{homotopy groups} of $G\hskip-.015in .$ are defined as
$$
\gathered
\pi_n(G\hskip-.015in .) = G^n_n/(\partial_{n+1}(G^n_{n+1}))\qquad n\ge 1\\
\pi_0(G\hskip-.015in .) = G_0/G^0_0
\endgathered
$$
These homotopy groups agree with the usual ones, in the sense that there is a natural
isomorphism of graded groups $\pi_*(G\hskip-.015in .) \cong \pi_*(|G\hskip-.015in .|,*)$
where the groups on the right are the homotopy groups of the base-pointed topological
space $|G\hskip-.015in .| =$ the geometric realization of $G\hskip-.015in .$, with basepoint
corresponding to $1\in G_0$. We call $G\hskip-.015in .$ a \underbar{resolution} if
$\pi_n(G\hskip-.015in .) = 0$ for $n > 0$. Similarly, if $(G\hskip-.015in .)_+$ is an
augmented simplicial group, it is a resolution if $G\hskip-.015in .$ is one, and
$\pi_0(G\hskip-.015in .) = G_{-1}$.

Throughout this paper, a $\underline{C^*-algebra}$ $A$ will always mean a normed,
involutive Banach algebra over $\Bbb C$ satisfying $\|x\|^2 = \|x^*x\|$. $C^*$-algebras
need not be unital. We denote by
$(C^*-algebras)$ the category of $C^*$-algebras and $C^*$-algebra homomorphisms.
A \underbar{simplicial $C^*$-algebra} is then a covariant functor $F:\Delta^{op}\to
(C^*-algebras)$, while an \underbar{augmented simplicial $C^*$-algebra} is a covariant
functor $F:\Delta^{op}_+\to (C^*-algebras)$. When referring to the homotopy groups of a
simplicial $C^*$-algebra $A\hskip-.015in .$, we will always mean its combinatorial
homotopy groups as defined above. Degreewise, this means ignoring the topology and
algebra structure on $A_n$, treating it as a discrete abelian group under addition.

If $A$ is a $C^*$-algebra, then $A^2 = A$ (cf. [Di, (1.5.8)]). In particular, if $I$ is a
closed $C^*$-ideal in a $C^*$-algebra $A$, it is a sub-$C^*$-algebra of $A$, so $I^2 = I$.

\proclaim{\bf\underbar{Proposition 1.3}} Every simplicial $C^*$-algebra $A\hskip-.015in .$
is a resolution.
\endproclaim

\prf For $n\ge 1$, $A^n_n$ is a closed $C^*$-ideal in $A_n$. However, for any simplicial
algebra $B\hskip-.015in .$, $(B^n_n)^2\subset \partial_{n+1}(B^n_{n+1})$ for all $n\ge 1$. To see this, suppose $a, b\in B^n_n$. Then $s_n(a)(s_n(b) - s_{n-1}(b))\in B^n_{n+1}$ and maps to $ab$ under $\partial_{n+1}$. Consequently
$\pi_n(A\hskip-.015in .) = 0$ for all $n \ge 1$.\hfill //
\endpf

Note also that as $A^0_0$ is the image of the closed $C^*$-ideal $A^0_1$ under $\partial_1$,
it is an involutive ideal in $A_0$ (by a simplicial argument), and closed (as it is the
image of a $C^*$-homomorphism). Hence the quotient $\pi_0(A\hskip-.015in .)$ is a quotient
$C^*$-algebra of $A_0$. Thus if $A\hskip-.015in .$ is a simplicial $C^*$-algebra, we may
define its associated augmented simplicial $C^*$-algebra $A\hskip-.015in .^+$ to be
$A\hskip-.015in .$ in non-negative degrees with $A_{-1}^+ = \pi_0(A\hskip-.015in .)$. By the
above proposition, any associated augmented simplicial $C^*$-algebra is a resolution.

By (1.3) we have short-exact sequences of $C^*$-algebras and ideals:
$$
\gathered
A^0_0\inj A_0\surj \pi_0(A\hskip-.015in .)\\
A^n_n\inj A^{n-1}_n\surj A^{n-1}_{n-1}
\endgathered
\tag1.4
$$

For a Banach algebra $B$ over $\Bbb C$, let $\underline{\underline{K}}(B)$ denote the
topological $K$-theory spectrum of $B$ [K]. We may define this as the $\Omega$ spectrum which
in even dimensions is $K_0(B)\times BGL(B)$ and in odd dimensions $GL(B)$. Then
$K_*^t(B) = \pi_*(\underline{\underline{K}}(B))$. The functor $B\mapsto
\underline{\underline{K}}(B)$ satisfies excision, in the sense that it associates to any
short-exact sequence of Banach algebras $I\inj B\surj \ov B$ a homotopy-fibration sequence
of spectra $\underline{\underline{K}}(I)\to \underline{\underline{K}}(B)\to
\underline{\underline{K}}(\ov B)$ (loc. cit.). Thus, given a simplicial $C^*$-algebra
$A\hskip-.015in .$, we may define a filtration on $K_*^t(\pi_0(A\hskip-.015in .))$ by
$$
\gathered
\Cal F_nK_*^t(\pi_0(A\hskip-.015in .)) =
ker(\partial^{(n)}: K_*^t(\pi_0(A\hskip-.015in .))\to K_{*-n}^t((A\hskip-.015in .)^{n-1}_{n-1})) = \\
ker(K_*^t(\pi_0(A\hskip-.015in .))\overset\partial\to\longrightarrow
K_{*-1}^t((A\hskip-.015in .)^0_0)\overset\partial\to\longrightarrow
K_{*-2}^t((A\hskip-.015in .)^1_1)\overset\partial\to\longrightarrow\dots
\overset\partial\to\longrightarrow K_{*-n}^t((A\hskip-.015in .)^{n-1}_{n-1}))
\endgathered
\tag1.5
$$

Recall that for a discrete group $\pi$, $C^*(\pi) = C^*_{max}(\pi)$ is defined as the
maximal $C^*$-completion of $\ell^1(\pi)$. The association $\pi\mapsto C^*(\pi)$
defines a functor from $(groups)$ (the category of discrete groups) to $(C^*-algebras)$.
Thus if $\Gamma\hskip-.03in .$ is a simplicial group, applying $C^*(_-)$ degreewise
produces the simplicial $C^*$-algebra  $C^*(\Gamma\hskip-.03in .) = \{C^*(\Gamma_n)\}_{n\ge 0}$.
Given $\Gamma\hskip-.03in .$, let $\pi = \pi_0(\Gamma\hskip-.03in .)$. Also associated to
$\Gamma\hskip-.03in .$ is the simplicial group algebra
$\Bbb C[\Gamma\hskip-.03in .]$. We will need the following result, due to Baum and Connes.

\proclaim{\bf\underbar{Lemma 1.6}} (Baum-Connes) Let $\phi : G_1\surj G_2$ be a surjective
homomorphism of discrete groups.
Let $I = \ker (\Bbb C[G_1]\overset \phi\to\surj\Bbb C[G_2])$, and let $I^*$ be the norm
closure of $I$ in
$C^*(G_1)$ under the natural embedding $\Bbb C[G_1]\inj C^*(G_1)$. Then there exists a
canonical isomorphism
of $C^*$ algebras $A = C^*(G_1)/I^* \overset\cong\to\to C^*(G_2)$.
\endproclaim
\prf First, $A$ is a $C^*$-algebra, as it is the quotient of a $C^*$-algebra by a norm-closed
and star-closed ideal. Because the image of any $C^*$-algebra homomorphism is closed [Di], it is surjective if and only if it has dense image. Thus the $C^*$-algebra homomorphism $C^*(\phi): C^*(G_1)\to C^*(G_2)$ induced by $\phi$ is surjective, as the image contains the dense subalgebra $\Bbb C[G_2]$. Also, $C^*(\phi)$ sends
$I^*$ to $0$, hence factors by $A$. Now
$(I^*\cap\Bbb C[G_1])\subseteq \ker(\Bbb C[G_1]\surj \Bbb C[G_2]\inj C^*(G_1)) = I$, so the
inclusion
$\Bbb C[G_1]\inj C^*(G_1)$ induces an inclusion
$\Bbb C[G_2] = \Bbb C[G_1]/(I^*\cap\Bbb C[G_1])\inj A$ which has dense image in the norm
topology.
Let $\rho : G_2\to \Cal B(\Cal H)$ be a unitary representation of $G_2$ on a complex
Hilbert space $\Cal H$.
Then $\phi\circ\rho : G_1\to \Cal B(\Cal H)$ is a unitary representation of $G_1$ which
by the universal property of the maximal $C^*$-algebra functor
$C^*(_-)$ admits a unique extension to a $*$-representation $C^*(G_1)\to \Cal B(\Cal H)$;
this representation sends $I^*$ to $0$, inducing a representation $A\to \Cal B(\Cal H)$.
Therefore $A$ satisfies the axioms which uniquely characterize $C^*(G_2)$, implying that
the $C^*$-algebra surjection $A\surj C^*(G_2)$ is an isomorphism.\hfill //
\endpf

Applying this lemma in the case $G_1 = \Gamma_0$ and $G_2 = \pi$ we conclude

\proclaim{\bf\underbar{Corollary 1.7}} The $C^*$-ideal
$C^*(\Gamma\hskip-.03in .)^0_0 \overset{def}\to = \partial_1(C^*(\Gamma\hskip-.03in .)^0_1)$ is the norm-closure of $\Bbb C[\Gamma\hskip-.03in .]^0_0 = ker(\Bbb C[\Gamma_0]\surj \Bbb C[\pi])$ in
$C^*(\Gamma_0)$. Consequently, $C^*(\Gamma\hskip-.03in .)$ is a resolution of
$\pi_0(C^*(\Gamma\hskip-.03in .)) = C^*(\pi)$ for any simplicial group
$\Gamma\hskip-.03in .$ with $\pi = \pi_0(\Gamma\hskip-.03in .)$.
\endproclaim

\prf By the previous Lemma, $C^*(\Gamma\hskip-.03in .)^0_1$ is the norm-closure of $\Bbb C[\Gamma\hskip-.03in .]^0_1$ in $C^*(\Gamma_1)$. Because the image of $C^*(\Gamma\hskip-.03in .)^0_1$ under $\partial_1$ is closed in $C^*(\Gamma_0)$, $C^*(\Gamma\hskip-.03in .)^0_0$ identifies with the norm-closure of $\partial_1(\Bbb C[\Gamma\hskip-.03in .]^0_1) = \Bbb C[\Gamma\hskip-.03in .]^0_0$ in $C^*(\Gamma_0)$.\hfill //
\endpf
\vskip.3in


\head 2 Filtering homotopy groups \endhead \vskip.2in

As noted above, an augmented simplicial object in a category $\underline C$ is a
covariant functor $F : \Delta^{op}_+ \to \underline C$. Let $n\ge-1, -1\le j\le n$.

\definition{\bf\underbar{Definition 2.1}} $\underline D^j_n \subset
\Delta^{op}_+$ is the subcategory with objects $\{\underline
n,\underline{n-1},\dots, \underline{n-j-1}\}$. The morphisms of $\underline
D^j_n$ are generated by
$$
\{\partial_i :
\underline{n-k} \to \underline{n-k-1}\; | \; 0\le i\le j-k\}\text{ for
}0\le k\le j,\text{ together with }Id_{\underline{n-k}}\; .
$$
\enddef Given a functor $F :
\Delta^{op}_+ \to \underline C$, it restricts to define a functor $F :
\underline D^j_n \to \underline C$. We will need the degeneracy morphisms
of $\Delta^{op}_+$ later on in verifying certain properties of this
restriction, however they are not needed for the following construction.

For an integer $m\ge 0$ , $C(\underline m)$ will denote the
category of subsets of $\{0,1,2,\dots, m\}$, with morphisms given
by inclusions. An $\underline{m+1 -\text{ cube }}$ in the category
$\underline C$ is by definition a covariant functor $F':
C(\underline m) \to \underline C$.

Let $F : \Delta^{op}_+ \to \underline C$ as before.

\proclaim{\bf\underbar{Lemma 2.2}} $F |_{\underline D^j_n} : \underline
D^j_n \to \underline C$ determines a $(j+1)-$cube in $\underline C$ :
$$
\gather F^j_n : C(\underline j) \to\underline C \;;\\
\noindent F^j_n \text{ is defined on morphisms and objects by}\\ F^j_n(S) =
F(\underline{n-|S|}), \;\phi \incl S \incl \underline j\;,\tag2.3\\
F^j_n(S\hookrightarrow S \cup \{i\}) = F(\partial_k) : F(\underline{n-|S|})
\to F(\underline{n-|S|-1})\;,
\endgather
$$
where $i \not\in S\, (0 \le i \le j)$, and $k= i-|S_i|$, with
$S_i = \{x
\in S| x < i\}$ .
\endproclaim

\prf We need to show
that all possible squares commute. Let $S \subseteq \underline j$ be fixed,
and assume $|S|\le j-1$. Choose $i_1, i_2$ with $0\le i_1,\; i_2 \le j,\;
i_1 \ne i_2$, and $S\cap
\{i_1, i_2\}=\phi$. Now consider the square
$$
\diagram S \rto|<<\tip \dto|<<\tip & S\cup \{i_1\} \dto|<<\tip \\ S \cup
\{i_2\} \rto|<<\tip & S\cup \{i_1, i_2\}
\enddiagram
\tagit{2.4}
$$

Applying $F^j_n$ we get a square
$$
\diagram F(\underline{n-|S|}) \rto^{\ssize{F(\partial_{k_1})}}
\dto_{\ssize{F(\partial_{k_2})}} & F(\underline{n-|S\cup\{i_1\}|})
\dto^{\ssize{F(\partial_{k_3})}} \\ F(\underline{n-|S\cup\{i_2\}|})
\rto^{\ssize{F(\partial_{k_4})}} & F(\underline{n-|S\cup\{i_1, i_2\}|})
\enddiagram
\tagit{2.5}
$$
where the subscripts $k_i$ are determined by
(2.3). We can assume without loss of generality that $i_1 < i_2$. Then
$$
\aligned k_1 &= i_1 - |S_{i_1}| = i_1 - |\{x\in S|x < i_1\}| \\
k_2 &= i_2 - |S_{i_2}| = i_2 - |\{x\in S|x < i_2\}| \\
k_3 &= i_2 - |(S\cup\{i_1\})_{i_2}| = i_2 - |\{x\in S\cup \{i_1\}| x < i_2\}| \\
k_4 &= i_1 - |(S\cup\{i_2\})_{i_1}| = i_1 - |\{x\in S\cup\{i_2\}| x < i_1\}|\;.
\endaligned
$$
The inequality $i_1 < i_2$ implies that $k_1 = k_4$ and $k_3
= k_2 - 1$. It also implies $k_3 \ge k_1$. The identity $\partial_{k_3}
\partial_{k_1} = \partial_{k_4}\partial_{k_2} $ verifies that the square
commutes. \hfill //
\enddemo

In what follows, $\underline C = \text{ (spaces)}_*, \text{(spectra)}_* \text{ or
the category (complexes)}_*$ of chain complexes over some fixed field. We define
\vskip.1in

\noindent(2.6)  i) $hp F^j_n =$ the iterated homotopy cofibre of
the $(j+1)-$ cube $F^j_n : C(\underline j) \to \underline C$
\vskip.02in

\indent\indent\; ii) $hf F^j_n =$ the iterated homotopy fibre of the
$(j+1)-$ cube $F^j_n$.
\vskip.1in

\noindent When the target is the category (spectra)$_*$ or (complexes)$_*$,  the natural
equivalence between
homotopy-fibration sequences and homotopy-cofibration sequences
produces an equivalence
$$
hf F^j_n \overset\simeq\to\longrightarrow \Omega^{j+1}\; hp\, F^j_n\; .
\tagit{2.7}
$$

\proclaim{\bf\underbar{Lemma 2.8}} Let $F:\Delta^{op}_+ \to \underline C$
. Then for all $n\ge 0$ and $-1\le j\le n-1$, there is a
homotopy-fibration sequence
$$ hf F^{j+1}_n
\overset{\overline\alpha^j_n}\to\longrightarrow hf F^j_n
\overset{\alpha^j_n}\to\longrightarrow (hf F^j_{n-1})
\tagit{2.9}
$$
where $\overline\alpha^j_n$ is defined below; moreover
there is an equivalence of maps
$$
\alpha^j_n \simeq \widetilde\partial_{j+1}: hf F^j_n \to hf F^j_{n-1}\;,
\tagit{2.10}
$$
where $\widetilde\partial_{j+1}$ is the map of iterated
homotopy fibres induced by $\partial_{j+1}$.
\endproclaim

\prf We begin by defining inclusions of categories

\noindent(2.11)\;\;\;\; i)\;\; $C(\underline j)
\overset{\tau_1}\to\longrightarrow C(\underline{j+1})$ ;

\indent\indent\quad $\tau_1(S) = S$ ,
\quad $\tau_1(S\hookrightarrow T) = (S\hookrightarrow T)$ .

\indent\indent\quad\phantom{1.} ii)\;\; $C(\underline j)
\overset{\tau_2}\to\longrightarrow C(\underline{j+1})$ ;

\indent\indent\quad $\tau_2(S) = S\cup \{j+1\}$ ,
\quad $\tau_2(S\hookrightarrow T) = (S\cup \{j+1\} \hookrightarrow T \cup
\{j+1\})$ .

 With $F$ as above, and $j<n$, define
$$
(F_i)^j_n : C(\underline j) \to
\underline C
\tagit{2.12}
$$
as the composition
$$
C(\underline j)
\overset{\tau_i}\to\longrightarrow C(\underline{j+1})
\overset{F^{j+1}_n}\to\longrightarrow \underline C \; .
$$

It is easy to see from (2.3) that $(F_1)^j_n = F^j_n \;,\; (F_2)^j_n =
F^j_{n-1}$. We may view the inclusion $\tau_2$ as a natural transformation from the
$(j+1)-$cube $(F_1)^j_n$ to the $(j+1)-$cube $(F_2)^j_n$; doing so yields a
homotopy-fibration sequence
$$
hfF^{j+1}_n \longrightarrow hf(F_1)^j_n
\longrightarrow hf(F_2)^j_n
\tagit{2.13}
$$
which, upon replacing $(F_1)^j_n$ by $F^j_n$ and
$(F_2)^j_n$ by $F^j_{n-1}$, agrees with the sequence in (2.9) and defines the
maps $\alpha^j_n$ , $\overline\alpha^j_n$. It remains to identify
$\alpha^j_n$ with $\widetilde\partial_{j+1}$ up to homotopy. As we are working up to
homotopy, we may assume without loss of generality that the $(j+2)$-cube has been
replaced in a functorial way by a $(j+2)$-cube in which all of the morphisms are mapped
to fibrations (and the homotopy fibres are fibres). In this case, there are natural
inclusions $hf(F_1)^j_n\hookrightarrow F(\phi)$, $hf(F_2)^j_n\hookrightarrow F(\{j+1\})$.
By (2.3)
$$
F(\phi\hookrightarrow \{j+1\}) = F(\partial_{j+1}) : F(\underline n) \to
F(\underline{n-1})\;,
\tagit{2.14}
$$
which identifies $\alpha^j_n$ up to homotopy as the map induced by $\partial_{j+1} $.\hfill
//
\endpf

\remark{\bf\underbar{Remark 2.15}} When $\underline C = \text{
(spaces)}_*$, care should be taken in interpreting the statement of the
Lemma 2.8. In our applications of this lemma, $F(\underline n)$ is a
connected space for each $n \ge -1$ . This implies $hfF^j_n$ is connected
for $j<n$ . For the homotopy fibration sequence $hfF^{j+1}_n
\longrightarrow hfF^j_n \longrightarrow hfF^j_{n-1}$ of (1.2.9)
yields a long-exact sequence in homotopy which terminates with
$\pi_1(hfF^j_{n-1})$ when $j<n-1$. However for $j=n-1, \;n>0$, the base space may
not be path-connected. In this case, the sequence should be written as
$$
hfF^n_n \longrightarrow hfF^{n-1}_n \longrightarrow (hfF^{n-1}_{n-1})_0 \; .
$$
If $\underline C = \text{(spectra)}_*$ or (complexes)$_*$, this
problem does not arise. Moreover in this case the natural equivalence
between homotopy-fibration and homotopy-cofibration sequences in this
category yield a dual formulation of (2.9) as a cofibration sequence
$$
hpF^j_n \longrightarrow hpF^j_{n-1} \longrightarrow hpF^{j+1}_n
\tagit{2.9${^*}$}
$$
where the first map is induced by $\partial_{j+1}$,
as before.

Now from the identification $F(\underline{-1}) \simeq hfF^{-1}_{-1}$,
we may use the sequences of (2.9) to define a filtration on
$\pi_*(F(\underline{-1}))$.
\endremark

\definition{\bf\underbar{Definition 2.16}} Let $F:\Delta^{op}_+ \to
\underline C$. For $k\ge 1$ set
$$
\Cal F_k P_* (F(\underline{-1}) ) =
\ker \left(\partial^{(k)}_* : P_* (hf\, F^{-1}_{-1})
\overset\partial\to\longrightarrow P_{*-1} (hf\, F^0_0)
\overset\partial\to\longrightarrow \dots
\overset\partial\to\longrightarrow P_{*-k} (hf\, F^{k-1}_{k-1})\right)
$$
where $P_*(_-) = \pi_*(_-)$ for  $\underline C = \text{(spaces)}_*$
or $\text{ (spectra)}_*$, and $H_*(_-)$ for  $\underline C = \text{(complexes)}_*$.
The maps which appear in the sequence are the boundary maps associated with
the sequence (2.9) for $-1 \le j=n-1 \le k-1$.
\enddef

This is clearly an increasing filtration of $\pi_* (F(\underline{-1}))$.
When $\underline C$ is the category of spaces, $\pi_* (F(\underline{- 1}))
= \underset k\to\varinjlim\,
\Cal F_k \pi_* (F(\underline{-1}))$ for dimensional reasons. When
$\underline C = \text{ (spectra)}_*$ or $\text{(complexes)}_*$,
the filtration of $\pi_* (F(\underline{-1}))$ is more natural
for the reason that long-exact sequences
are allowed to continue on indefinitely, preserving exactness.
Because of the functoriality of the constructions involved, we
have

\proclaim{\bf\underbar{Proposition 2.17}} A natural transformation $F_1
\overset\zeta\to\longrightarrow F_2:
\Delta ^{op}_+ \to \underline C$ induces compatible natural transformations
$(F_1)^j_n \to (F_2)^j_n$, and hence a filtration-preserving homomorphism
of homotopy groups
$$
\Cal F_*\zeta_* : \Cal F_* \pi_* (F_1(\underline{-1})) \longrightarrow
\Cal F_* \pi_* (F_2(\underline{-1}))\; .
$$
\endproclaim

\prf Clear. \hfill //
\endpf

Finally, we note the following fact, which will be used often in what follows.

\proclaim{\bf\underbar{Proposition 2.18}} If
$F : \Delta ^{op}_+ \to\underline C =$ (complexes)$_*$,
and for each $k$ the simplicial abelian group $\{[n]\mapsto ((F(n))_k\}_{n\ge 0}$
is a resolution (i.e., has vanishing homotopy groups above dimension $0$), then for
each $-1\le j\le n$ there is a weak equivalence
$$
fF^j_n\overset\simeq\to\inj hfF^j_n
$$
where $fF^j_n$ denotes the iterated fibre of the $(j+1)$- cube $F^j_n$, and the
map is the natural inclusion of the iterated fibre into the iterated homotopy
fibre.
\endproclaim

\prf This follows by induction on $j$ starting with $j= -1$ and showing that for
each $k$
$$
fG^{j+1}_n
\overset{\overline\alpha^j_n}\to\longrightarrow fG^j_n
\overset{\alpha^j_n}\to\longrightarrow fG^j_{n-1}
$$
is a fibration sequence, where $G(n) = F(n)_k$.
For $j<n$ this follows immediately from the simplicial
identities which show that the map on the right is a surjection of simplicial
abelian groups. In the final case $j=n$, the surjectivity of this map follows
from the hypothesis on $F$.\hfill //
\endpf
\vskip.3in


\head 3 A basic example
\endhead
\vskip.2in

We will identify the filtration of the preceding section in a simple but
important case.

For a simplicial set $X\hskip-.015in .$, let $\BbbZ\{X\hskip-.015in .\}$ denote the free
simplicial abelian group generated by $X\hskip-.015in .$. The chain complex associated to
$\BbbZ\{X\hskip-.015in .\}$ is exactly the singular chain complex for $X\hskip-.015in .$,
so $\pi_*(\BbbZ\{X\hskip-.015in .\})\cong H_*(X\hskip-.015in .)$. For a discrete group $G$,
$BG$ will denote the standard simplicial (non-homogeneous) bar construction on $G$.

In this section $\Gamma\hskip-.03in .$ denotes a free simplicial resolution of
a discrete group $\pi$.  Define $G^j_nB(\gm\hskip-.03in .)$ for $n\ge 1$ by
$$
G^j_nB(\gm\hskip-.03in .) = \intjn \BbbZ\{B\gm_n\}
\to \BbbZ\{B\gm_{n-1}\})\; .
\tag3.1
$$
Thus $G^j_nB(\gm\hskip-.03in .)$ is a simplicial abelian subgroup of
$\BbbZ\{B\gm_n\}$ for $0\le j\le n$ , $n\ge 1$. Similarly we define
$$
\gather G^0_0B(\gm\hskip-.03in .) = \ker (\BbbZ\{B\gm_0\} \surj \BbbZ\{B\pi\})\;,\\
G^{-1}_nB(\gm\hskip-.03in .) = \BbbZ\{B\gm_n\}\;\text{ for }n\ge 0\;,\tag3.2\\
G^{-1}_{-1}B(\gm\hskip-.03in .) = \BbbZ\{B\pi\}\; .
\endgather
$$

\proclaim{\bf\underbar{Lemma 3.3}} i) For each $n\ge 0$ and $-1
\le j\le n-1$, there is a short-exact sequence of simplicial abelian groups
$$
G^{j+1}_nB(\gm\hskip-.03in .) \rightarrowtail G^j_nB(\gm\hskip-.03in .)
\overset{(\partial_{j+1})_*}\to\twoheadrightarrow G^j_{n-1}B(\gm\hskip-.03in .)
\tag3.4
$$

\indent\indent\indent ii)\quad For $0\le j\le n,\; \pi_0(G^j_nB(\gm\hskip-.03in .)) = 0$ .
For
$0\le j\le n-1, \; \pi_i(G^j_nB(\gm\hskip-.03in .))=0$ when $i>1$.
\endproclaim

\prf Certainly $G^{j+1}_nB(\gm\hskip-.03in .)$ is the kernel of $(\partial_{j+1})_*$
restricted to $G^j_nB(\gm\hskip-.03in .)$. The simplicial identities
imply that $ (\partial_i)_*\circ(\partial_{j+1})_* = (\partial_j)_* \circ
(\partial_i)_* $ for $i \le j$, so the image of
$(\partial_{j+1})_*\,|_{G^j_nB(\gm\hskip-.03in .)}$ is contained in
$G^j_{n-1}B(\gm\hskip-.03in .)$.

If $j < n - 1$ and $x \in G^j_{n-1}B(\gm\hskip-.03in .)$, then $(s_{j+1})_*(x) \in
G^j_nB(\gm)$ and maps to $x$ under $(\partial_{j+1})_*$, proving
surjectivity in this case.

Now fix $p$, and consider the simplicial abelian group $\{\BbbZ
\{(B\gm_n)_p\}\}_{n \ge 0}$ . $(B\gm_n)_p \cong (\gm_n)^p$, and the
simplicial set $\{(B\gm_n)_p\}_{n\ge 0}$ is isomorphic as a simplicial set
to the simplicial group $\{(\gm_n)^p\}_{n\ge 0} \cong (\gm\hskip-.03in .)^p$ equipped
with diagonal simplicial structure. Thus the
simplicial abelian group $\{\BbbZ\{(B\gm_n)_p\}\}_{n\ge 0}$ identifies with
the underlying abelian group of the simplicial group algebra $\BbbZ[(\gm\hskip-.03in .)^p]$.
As
$\gm\hskip-.03in . \twoheadrightarrow \pi$ is a weak equivalence, $(\gm\hskip-.03in .)^p
\overset\cong\to\longrightarrow \pi^p$ is again a weak equivalence for each $p>0$.
Thus $\pi_i(\BbbZ[(\gm\hskip-.03in .)^p]) = 0$ for all $i > 0$ and
$p>0$. The first part of the lemma now follows by the same argument as in
Proposition 2.18.

For each $n\ge 0$, $B\gm_n$ is a $0$-reduced simplicial set. Thus
$\BbbZ\{B\gm_n\}_0 \cong \BbbZ$ for all $n\ge 0$; moreover the simplicial
structure on $\{\BbbZ\{B\gm_n\}_0\}_{n\ge0}$ induced by the face and
degeneracy maps of $\gm\hskip-.03in .$ is the trivial one. It follows that
$G^j_nB(\gm\hskip-.03in .)_0
= 0$ for all $n, \; 0\le j \le n$ so $\pi_0(G^j_nB(\gm\hskip-.03in .)) = 0$ . Now
suppose $n\ge 1$.

The sequence
$$
G^0_nB(\gm\hskip-.03in .) \longrightarrow
\BbbZ\{B\gm_n\}\overset{(\partial_0)_*}\to\twoheadrightarrow
\BbbZ\{B\gm_{n-1}\}
\tag3.5
$$
is split-exact $((s_0)_*$ splits $(\partial_0)_*)$.
$\gm_j$ is a free group for each $j$  so $\pi_i(\BbbZ\{B\gm_j\}) = H_i(B\gm_j) =0$ for
$i>1,\; j\ge 0$. The split-exactness of (3.5) implies
$\pi_i(G^0_nB(\gm\hskip-.03in .))=0$ for $i>1,\; n\ge 0$. Inductively, assume the
result for $j-1$. By the above there is a short exact
sequence
$$
G^j_nB(\gm\hskip-.03in .)
\rightarrowtail G^{j-1}_n B(\gm\hskip-.03in .)
\overset{(\partial_j)_*}\to\twoheadrightarrow G^{j-1}_{n-1}B(\gm\hskip-.03in .)\; .
\tag3.6
$$
which is split-exact for $j<n$. By the five lemma, $\pi_i(G^j_nB(\gm\hskip-.03in .)) = 0$
for $i>1$.
\hfill //
\endpf

We define functors $F_i, \; i = 0,1,2$ by

\noindent (3.7)
\roster
\item $F_0: \Delta ^{op}_+ \to (groups)$, $F_0(\underline n) = \gm_n
\;\; n\ge 0$ , $F_0 (\underline{-1}) = \pi$ . Thus $F_0|_{\Delta
^{op}}$ is the simplicial group $\gm\hskip-.03in .$ and $F_0 (\underline 0 \to
\underline{-1}) = \epsilon: \gm_0\twoheadrightarrow \pi$.

\item $F_1 = \BbbZ \{BF_0\}$. Thus for each $n\ge-1,\; F_1(\underline n)$
is the simplicial abelian group $\BbbZ\{BF_0 (\underline n)\}$.

\item $F_2$ is the augmented simplicial functor which assigns to $n$ the
associated chain complex of $F_1(\underline n)$, with face and degeneracy maps induced
by $\gm\hskip-.03in .$.
\endroster

By the previous section, the restriction of $F_i$ to the subcategory
$\underline D^j_n \subset \Delta ^{\text{op}}_+$ determines a $(j+1)-$cube
$(F_i)^j_n : C(\underline j) \to \underline C$, where $\underline C$ is the
target of $F_i$.

\proclaim{\bf\underbar{Lemma 3.8}} For all $n \ge -1,\; -1 \le j\le n$ ,
there are weak equivalences
$$
\gather
hf(F_1)^j_n \overset\simeq\to\longleftarrow f(F_1)^j_n = G^j_nB(\gm\hskip-.03in .)\\
hf(F_2)^j_n \overset\simeq\to\longleftarrow f(F_2)^j_n = (G^j_nB(\gm\hskip-.03in .))_*
\endgather
$$
where $(A\hskip-.015in .)_*$ denotes the associated chain complex of the simplicial
abelian group $A\hskip-.015in .$. Under this equivalence, the homotopy
fibration sequence of (2.9) is identified with the fibration sequence of
(3.4), up to homotopy.
\endproclaim

\prf By definition, $hf(F_1)^{-1}_n\overset =\to\longrightarrow G^{-1}_nB(\gm\hskip-.03in .)$.
By Proposition 2.18, the hypothesis on $\gm\hskip-.03in .$ implies
that the natural map from the iterated fibre -- $(G^j_nB(\gm\hskip-.03in .))_*$ --
to the iterated homotopy fibre $hf(F_1)^j_n$ is a weak equivalence for all $j\le n$. The
same argument works for $F_1$.\hfill //
\enddemo

The construction in Definition 2.16 produces a filtration on
$\pi_*(F_1(\underline{-1})) = H_*(B\pi)$ which may alternatively be defined
using the fibration sequences in (3.4).
Lemma 3.3 and the fact that
$\pi_i(G^j_nB(\gm\hskip-.03in .)) =0$ for $i > 0$ and $j<n$ implies that\newline
$\partial:\pi_i(G^{k-1}_{k-1}B(\gm\hskip-.03in .)) \rightarrow
\pi_{i-1} (G^k_kB(\gm\hskip-.03in .))$
is an isomorphism if $i>2$ and injective if $i = 2$. Thus $\pi_j(\BbbZ\{B\pi\})$
maps injectively by $\partial^{(j-1)}$ to $\pi_1(G^{j-2}_{j-2}B(\gm\hskip-.03in .))$,
and then
to zero under
$\partial:\pi_1(G^{j-2}_{j-2}B(\gm\hskip-.03in .)) \rightarrow
\pi_0(G^{j-1}_{j-1}B(\gm\hskip-.03in .)) = 0$.
It follows that $0 \ne x \in \pi_j (\BbbZ\{B\pi\})$ maps to zero under
$\partial^{(n)}$ if and only if $j \le n$. The same remarks apply to the functor
$F_2$. We conclude

\proclaim{\bf\underbar{Proposition 3.9}} The filtration on
$$
H_*(F_2(\underline{-1})) = \pi_*(F_1(\underline{-1})) = H_*(B\pi;\Bbb Z)
$$
defined by (2.16) agrees with that induced by the skeletal filtration of $B\pi$:
$$
\Cal F_n H_*(F_2(\underline{-1})) =
\Cal F_n \pi_*(F_1(\underline{-1})) = \Cal F_n H_* (B\pi;\Bbb Z) =
\bigoplus_{i=0}^n H_i (B\pi;\Bbb Z)\; .
$$
\endproclaim
\vskip.3in


\head 4 Filtering the assembly map
\endhead
\vskip.2in

As above, let $\Gamma\hskip-.03in .^+$ be an augmented free simplicial resolution of $\pi$.
Define functors $F_i : \Delta_+^{op}\to (spectra)_*$ by $F_1(\underline n) =
\underline{\underline{KU}}(B\Gamma_n) =$ the (unreduced) 2-periodic complex $K$-homology spectrum
of the classifying space $B\Gamma_n$, and
$F_2(\underline n) = \underline{\underline{K}}^t(C^*(\Gamma_n)) =$ the
2-periodic topological $K$-theory spectrum of $C^*(\Gamma_n)$. The assembly map determines a
natural transformation
$$
\Cal A :F_1\to F_2
\tag4.1
$$
By (2.16), $\Gamma\hskip-.03in .^+$ determines filtrations on both
$\pi_*(F_1(\underline{-1})) = KU_*(B\pi)$ and $\pi_*(F_2(\underline{-1})) =
K_*^t(C^*(\pi))$, and by Proposition 2.17, $\Cal A_*(\pi) =$ the map on homotopy groups
induced by $\Cal A(\pi)$, is filtration-preserving.

By (1.5) and (1.7), the simplicial $C^*$-algebra $C^*(\Gamma\hskip-.03in .^+)$ also
induces a filtration on $K_*^t(C^*(\pi))$. Our first observation is that that this
filtration agrees with that given by (2.16). Precisely, it follows from excision in
topological $K$-theory that for each $-1\le k\le n$ there is a commuting diagram
$$
\diagram
\underline{\underline{K}}^t(C^*(\Gamma\hskip-.03in .^+)^k_n)\dto^{\simeq}\rto &
\underline{\underline{K}}^t(C^*(\Gamma\hskip-.03in .^+)^{k-1}_n)\dto^{\simeq}\rto^{\wt{\partial_k}} &
\underline{\underline{K}}^t(C^*(\Gamma\hskip-.03in .^+)^{k-1}_{n-1})\dto^{\simeq}\\
hf(F_2)^k_n\rto & hf(F_2)^{k-1}_n\rto^{\wt{\partial_k}} &  hf(F_2)^{k-1}_{n-1}
\enddiagram
$$
where each row is a homotopy-fibration sequence and each vertical arrow is a weak equivalence.
This implies the equality of the sequence in (1.5) with that in (2.16). The assembly
map induces maps of homotopy-fibers
$$
hf\Cal A^k_n : hf(F_1)^k_n\to hf(F_2)^k_n
\tag4.2
$$
For $k = -1$ and $n\ge 0$ this is the assembly map for a free group, which is a weak
equivalence. The homotopy-fibration sequence of (2.9) is naturally split up to homotopy
by the degeneracy $s_{j+1}$. Thus $hf\Cal A^k_n$ is a weak equivalence for all $-1\le k < n$.
Together with an easy diagram-chase, this implies that if $x\in KU_*(B\pi)$ and
$\Cal A_*(\pi)(x)\in \Cal F_n K_*^t(C^*(\pi))$, then $x\in \Cal F_n KU_*(B\pi)$ (in fact
this follows simply from the surjectivity of $\Cal A_*(F)$ for free groups $F$). Finally,
via the functorial Atiyah-Hirzebruch Chern isomorphism
$KU_*(B\pi)\otimes\Bbb Q\overset{ch}\to{\underset\simeq\to\longrightarrow}
\bigoplus_{i\in\Bbb Z}H_{*-2i}(B\pi;\Bbb Q)$,
the computation in the previous section identifies, at least rationally, the filtration
on $KU_*(B\pi)$ induced by $\Gamma\hskip-.03in .^+$ as the skeletal filtration
$$
\Cal F_n KU_*(B\pi)\otimes\Bbb Q =
\underset{*-2i\le n}\to{\underset {i\in\Bbb Z}\to\bigoplus} H_{*-2i}(B\pi;\Bbb Q)
\tag4.3
$$

Summarizing, we have

\proclaim{\bf\underbar{Theorem 4.4}} Any free augmented simplicial resolution
$\Gamma\hskip-.03in .^+$ of a discrete group $\pi$ induces a filtration of both
$KU_*(B\pi)$ and $K_*^t(C^*(\pi))$. The assembly map preserves this filtration. The
filtration on $KU_*(B\pi)$ identifies rationally with the skeletal filtration (4.3),
while the filtration on $K_*^t(C^*(\pi))$ is defined as in (1.5). Finally, if
$x\in KU_*(B\pi)$ and $\Cal A_*(\pi)(x)\in \Cal F_n K_*^t(C^*(\pi))$, then
$x\in \Cal F_n KU_*(B\pi)$.
\endproclaim

As indicated in the introduction, the assembly map $\Cal A_*(\pi)$ is rationally injective
iff it is rationally injective on $ku_*(B\pi)\otimes\Bbb Q$, the rational connective
$K$-homology groups of $B\pi$. As a first corollary to the previous theorem we have

\proclaim{\bf\underbar{Corollary 4.5}} The assembly map resticted to connective
$K$-homology $\Cal A_*(\pi) : ku_*(B\pi)\to K_*^t(C^*(\pi))$ is rationally injective iff
the composition
$$
H_n(B\pi;\Bbb Q)\inj ku_n(B\pi)\otimes\Bbb Q\overset \Cal A_n(\pi)_{\Bbb Q}
\to\longrightarrow K_n^t(C^*(\pi))
\tag4.6
$$
is injective for all $n\ge 1$, where the first map is the inclusion of the top-dimensional
summand under the isomorphism $ku_n(B\pi)\otimes\Bbb Q\cong
\oplus_{i\ge 0}H_{n-2i}(B\pi;\Bbb Q)$.
\endproclaim

\prf We assume given an element of $ku_n(B\pi)\otimes\Bbb Q$, written as
$x = (x_n,x_{n-2},\dots,x_{\varepsilon})$ where $x_n\ne 0$, $x_m\in H_m(B\pi;\Bbb Q)$ and
$\varepsilon = 0$ or $1$. The image of this element under the rationalized assembly map
lies in $\Cal F_n K_n^t(C^*(\pi))\otimes\Bbb Q$, and under the composition
$$
\gather
ku_n(B\pi)\otimes\Bbb Q = \Cal F_n KU_n(B\pi)\otimes\Bbb Q\overset \Cal A_n(\pi)_{\Bbb Q}\to
\longrightarrow \Cal F_n K_n^t(C^*(\pi))\otimes\Bbb Q\\
\surj
\left( \Cal F_n K_n^t(C^*(\pi))\otimes\Bbb Q\right)/\left( \Cal F_{n-1} K_n^t(C^*(\pi))
\otimes\Bbb Q\right)
\endgather
$$
has the same image as the element $\wt x = (x_n,0,0,\dots,0)$. Suppose
$\Cal A_n(\pi)_{\Bbb Q}(\wt x)\ne 0$ in $\Cal F_n K_n^t(C^*(\pi))\otimes\Bbb Q$, but maps
to zero in the quotient
$\left( \Cal F_n K_n^t(C^*(\pi))\otimes\Bbb Q\right)/\left( \Cal F_{n-1} K_n^t(C^*(\pi))
\otimes\Bbb Q\right)$.
This would imply
$\Cal A_n(\pi)_{\Bbb Q}(\wt x)\in \Cal F_{n-1}  K_n^t(C^*(\pi))\otimes\Bbb Q$, which by the
previous theorem would imply $\wt x\in \Cal F_{n-1}  KU_n(C^*(\pi))\otimes\Bbb Q$ - a
contradiction (again by the previous theorem). Therefore if $\wt x$ maps non-trivially
to $\Cal F_n K_n^t(C^*(\pi))\otimes\Bbb Q$, it maps non-trivially to
$\left( \Cal F_n K_n(C^*(\pi))\otimes\Bbb Q\right)/\left( \Cal F_{n-1} K_n(C^*(\pi))
\otimes\Bbb Q\right)$, and therefore so does $x$. The other direction is obvious.\hfill //
\endpf

Referring to (1.5), we see that
$\left( \Cal F_n K_n^t(C^*(\pi))\otimes\Bbb Q\right)/\left( \Cal F_{n-1} K_n^t(C^*(\pi))
\otimes\Bbb Q\right)$ injects into $K_1^t(C^*(\Gamma\hskip-.03in .^+)^{n-2}_{n-2})
\otimes\Bbb Q$, and then maps to zero under the boundary map to
$K_0^t(C^*(\Gamma\hskip-.03in .^+)^{n-1}_{n-1})\otimes\Bbb Q$. In other words, it factors
by the inclusion
$$
coker\left(K_1^t(C^*(\Gamma\hskip-.03in .^+)^{n-1}_{n-1})\otimes\Bbb Q\to
K_1^t(C^*(\Gamma\hskip-.03in .^+)^{n-2}_{n-1})\otimes\Bbb Q\right)\inj
K_1^t(C^*(\Gamma\hskip-.03in .^+)^{n-2}_{n-2})\otimes\Bbb Q
\tag4.7
$$

Thus the last corollary can be further refined to

\proclaim{\bf\underbar{Corollary 4.8}} The composition
$$
\gathered
H_n(B\pi;\Bbb Q)\inj \Cal F_n KU_n(B\pi)\otimes\Bbb Q\overset \Cal A_n(\pi)_{\Bbb Q}
\to\longrightarrow \Cal F_n K_n^t(C^*(\pi))\\
\to
coker\left(K_1^t(C^*(\Gamma\hskip-.03in .^+)^{n-1}_{n-1})\otimes\Bbb Q
\to K_1^t(C^*(\Gamma\hskip-.03in .^+)^{n-2}_{n-1})\otimes\Bbb Q\right)
\endgathered
\tag4.9
$$
is injective for a given $n$ if and only if (4.6) is injective for that $n$.
\endproclaim

In terms of absolute $K$-groups, another corollary of the above is

\proclaim{\bf\underbar{Colollary 4.10}} The Strong Novikov Conjecture holds for $C^*(\pi)$ if and only if the composition
$$
H_n(B\pi;\Bbb Q)\overset \Cal A(\pi)_{\Bbb Q}\to\longrightarrow
K_n^t(C^*(\pi))\otimes\Bbb Q\overset\partial^{(n-1)}\to\longrightarrow
K_1^t(C^*(\pi)^{n-2}_{n-2})\otimes\Bbb Q
$$
is injective for each $n\ge 1$
\endproclaim

One possible method for detecting this image is via the map on $K_1^t(_-)$ induced by the surjection
$C^*(\pi)^{n-2}_{n-2}\surj C^*_{ab}(\pi)^{n-2}_{n-2}$ where $C^*_{ab}(\pi)^{n-2}_{n-2}$ denotes the $C^*$-algebra abelianization of $C^*(\pi)^{n-2}_{n-2}$ (i.e., the quotient by the norm-closure of the commutator ideal).
\vskip.3in


\head 5 Simplicial $C^{i}$-algebras and the Strong Novikov Conjecture
\endhead
\vskip.2in

 Let $G$ be a discrete group. We define $C^i(G)$ as the inverse-closure of $\Bbb C[G]$
 in $C^*(G)$. In other words, i) $C^i(G)$ contains $\Bbb C[G]$ and ii) if $a\in C^i(G)$
 and $a^{-1}\in C^*(G)$, then $a^{-1}\in C^i(G)$. This construction is clearly functorial
 with respect to group homomorphisms. By a [Sch], $M_n(C^i(G))$ is then
 inverse-closed in $M_n(C^*(G))$ for all $n$, where $M_n(A)$ denotes the ring of
 $n\times n$ matrices with coefficients in $A$. From the resulting equality
 $$
 GL_n(C^i(G)) = M_n(C^i(G))\cap GL_n(C^*(G))\qquad n\ge 1
 $$
 we conclude the inclusion $C^i(G)\hookrightarrow C^*(G)$ induces a weak equivalence
 of topological groups
 $$
 GL(C^i(G))\overset\simeq\to\hookrightarrow GL(C^*(G))
 \tag5.1
 $$
 and hence an isomorphism of higher topological $K$-groups
 $$
 K_*^t(C^i(G))\overset\cong\to\to K_*^t(C^*(G))\qquad *\ge 1
 $$
 (the topology on $C^i(G)$ is the fine topology, unless stated otherwise. Because $C^i(G)$ is
 inverse-closed in $C^*(G)$, its higher topological $K$-groups are the same for any topology
 between the fine topology and the induced norm topology [O1, App.]).
 The assignment $G\mapsto C^i(G)$ is functorial with respect to group homomorphisms,
 hence extends to simplicial and augmented simplicial groups.

 As in the previous section, $\Gamma\hskip-.03in .^+$ denotes a free augmented simplicial
 resolution of $\Gamma_{-1} = \pi$, and $C^i(\Gamma\hskip-.03in .^+)$ the resulting
 augmented simplicial algebra formed by applying $C^i(_-)$ degreewise. The corresponding
 simplicial group $\Gamma\hskip-.03in .$ is gotten by omitting the degree $-1$ part.

 \proclaim{\bf\underbar{Lemma 5.2}} There is a factorization
$$
\gathered
H_n(B\pi;\Bbb Q)\overset \Cal A_n(\pi)_{\Bbb Q}\to\longrightarrow \Cal F_n K_n^t(C^*(\pi))\\
\to
coker\left(K_1^t(\partial_n(C^i(\Gamma\hskip-.03in .^+)^{n-1}_{n}))\otimes\Bbb Q
\to K_1^t(C^i(\Gamma\hskip-.03in .^+)^{n-2}_{n-1})\otimes\Bbb Q\right)\\
\overset\cong\to\longrightarrow coker\left(K_1^t(C^*(\Gamma\hskip-.03in .^+)^{n-1}_{n-1})
\otimes\Bbb Q\to K_1^t(C^*(\Gamma\hskip-.03in .^+)^{n-2}_{n-1})\otimes\Bbb Q\right)
\endgathered
\tag5.3
$$
\endproclaim

\prf The inclusion $C^i(\Gamma\hskip-.03in .^+)\hookrightarrow C^*(\Gamma\hskip-.03in .^+)$
of augmented simplicial algebras yields inclusions
$C^i(\Gamma\hskip-.03in .^+)^k_n\hookrightarrow C^*(\Gamma\hskip-.03in .^+)^k_n$ which
induce isomorphisms in higher topological $K$-theory
$$
K_*(C^i(\Gamma\hskip-.03in .^+)^k_n)\cong K_*(C^*(\Gamma\hskip-.03in .^+)^k_n)\qquad *\ge 1
$$
The case $k < n$ follows from the case $k = -1$. When $k = n$, it is clear that
$C^i(\Gamma\hskip-.03in .^+)^n_n$ is inverse-closed in $C^*(\Gamma\hskip-.03in .^+)^n_n$.
By Lemma 1.6, it is dense in the norm topology. Now the inclusion
$\partial_n(C^i(\Gamma\hskip-.03in .^+)^{n-1}_{n})\hookrightarrow
C^i(\Gamma\hskip-.03in .^+)^{n-1}_{n-1}$ may not be surjective (in fact, this is an open
question; see below). However, it also induces an isomorphism in higher topological
$K$-theory. This is because  $\partial_n(C^i(\Gamma\hskip-.03in .^+)^{n-1}_{n})
\hookrightarrow \partial_n(C^*(\Gamma\hskip-.03in .^+)^{n-1}_{n})$ does, and
$\partial_n(C^*(\Gamma\hskip-.03in .^+)^{n-1}_{n}) = C^*(\Gamma\hskip-.03in .^+)^{n-1}_{n-1}$.
The result follows from (4.7).\hfill //
\endpf

In conjunction with section 1, the previous lemma implies

\proclaim{\bf\underbar{Theorem 5.4}} The Strong Novikov Conjecture is true for $\pi$
if and only if the map in (5.3)
$$
H_n(B\pi;\Bbb Q)\overset\wt{\Cal A}^t(\pi)_n\to\longrightarrow
coker\left(K_1^t(\partial_n(C^i(\Gamma\hskip-.03in .^+)^{n-1}_{n}))\otimes\Bbb Q
\to K_1^t(C^i(\Gamma\hskip-.03in .^+)^{n-2}_{n-1})\otimes\Bbb Q\right)
$$
is injective for each $n\ge 1$.
\endproclaim

We denote by $K(\Bbb Q,n-1)\hskip-.015in .$ the simplicial abelian
group which is $\{id\}$ in degrees less than $n-1$, equal to $\Bbb
Q$ in degree $n-1$, and equal to $\underset{s\in
S_{n-1,m}}\to\oplus \Bbb Q_s$ in degree $m\ge n$, where
$S_{n-1,m}$ is the set of distinct iterated degeneracies from
dimension $n-1$ to $m$, and $\Bbb Q_s$ indicates the copy of $\Bbb
Q$ gotton by applying $s$ to the unique copy of $\Bbb Q$ in
dimension $n-1$. There are isomorphisms
$$
H^n(B\pi;\Bbb Q) = [B\pi,K(\Bbb Q,n)]_* = [|B\Gamma\hskip-.03in .|,K(\Bbb Q,n)]_*
\cong Hom_{s.gps}(\Gamma\hskip-.03in ., K(\Bbb Q,n-1)\hskip-.015in .)
\tag5.5
$$
where $[X,Y]_*$ denotes (basepointed) homotopy classes of maps between (basepointed)
spaces $X$ and $Y$, while the right-hand side denotes the (abelian group) of simplicial
group homomorphisms from $\Gamma\hskip-.03in .$ to $K(\Bbb Q,n-1)\hskip-.015in .$.
The first two equalities are well-known, while the third follows by an easy simplicial
argument. And as we have noted, $C^i(_-)$ is functorial, so a simplicial group
homomorphism $\phi : \Gamma\hskip-.03in .\to K(\Bbb Q,n-1)\hskip-.015in .$ induces a
continuous homomorphism of $C^i$-algebras
$$
\wt\phi: C^i(\Gamma\hskip-.03in .)\to C^i(K(\Bbb Q,n-1)\hskip-.015in .)
\tag5.6
$$

The following useful fact was communicated to the author by P.
Baum [B].

\proclaim{\bf\underbar{Lemma 5.7}} For a discrete group $G$, let
$I_{\Bbb C}[G] = ker(\varepsilon:\Bbb C[G]\to \Bbb C)$, and $I^i(G) = ker(\varepsilon:C^i(G)
\to \Bbb C)$. Then the canonical inclusion $\Bbb C[G]\hookrightarrow C^i(G)$ induces an
isomorphism
$$
I_{\Bbb C}[G]/(I_{\Bbb C}[G])^2\overset\cong\to\longrightarrow I^i(G)/(I^i(G))^2
\tag5.8
$$
\endproclaim

\prf Because the map $G\to G_{ab} = G/[G,G]$ induces an isomorphism on the respective
quotients in (5.8), it suffices to consider the case when $G$ is abelian, and one may
further restrict to the case $G$ is torsion-free. By considering elements in $\Bbb C[G]$
of the form $1 + a$ where $a\in I_{\Bbb C}[G]$ and $|a| < 1$ ($|a| =$ the $C^*$-norm of $a$)
one sees the map in (5.8) is surjective. Injectivity follows from the finitely-generated
case by a direct limit argument. Finally, the finitely-generated torsion-free case reduces
to the case $G = \Bbb Z$ which is easy to verify directly.\hfill //
\endpf

Suppose given a homomorphism $\wt\phi$ of simplicial algebras as in (5.6). It is easy to
see that $C^i(K(\Bbb Q,n-1)\hskip-.015in .)_{n-1}^{n-1} = I^i(\Bbb Q)$ and
$\partial_n(C^i(K(\Bbb Q,n-1)\hskip-.015in .)^{n-1}_{n}) = (I^i(\Bbb Q))^2$. Let $K^a_*(R)$
denote the algebraic $K$-groups of the discrete ring $R$. We then have an induced
homomorphism
$$
\gathered
(\wt\phi)^a_1: coker\left(K_1^a(\partial_n(C^i(\Gamma\hskip-.03in .^+)^{n-1}_{n}))\otimes\Bbb Q
\to K_1^a(C^i(\Gamma\hskip-.03in .^+)^{n-2}_{n-1})\otimes\Bbb Q\right)\\
\longrightarrow
coker\left(K_1^a(\partial_n(C^i(K(\Bbb Q,n-1)\hskip-.015in .)^{n-1}_{n}))\otimes\Bbb Q
\to K_1^a(C^i(K(\Bbb Q,n-1)\hskip-.015in .)^{n-2}_{n-1})\otimes\Bbb Q\right)\\
\longrightarrow K_1^a(I^i(\Bbb Q)/(I^i(\Bbb Q))^2)\otimes\Bbb Q\cong\Bbb C
\endgathered
\tag5.9
$$

This composition is the analogue in algebraic $K$-theory which, if
it could be extended to topological $K$-theory would imply
injectivity of the composition in (1.5.4) in all dimensions. Our
next task is to formulate a version of the statement in Theorem
5.4 which is true in algebraic $K$-theory. To do so, we first give
an alternative description of $H_*(B\pi;\Bbb Q)$. Let
$\{D_n\}_{n\ge 0}$ be the chain complex $D_n =
H_1(B\Gamma_{n-1}^{n-2};\Bbb Q)$, $d^D_n = (\partial_{n-1})_* :
D_n\to D_{n-1}$. For $n<1$ set $D_n = 0$. For each $n \ge 1$,
$\Gamma_{n-1}$ acts on $\Gamma_{n-1}^{n-2}$ by conjugation. This
defines an action of $\Gamma_{n-1}$ on $D_n$; let $\ov{D}_n =
(D_n)_{\Gamma_{n-1}}$, the coinvariants with respect to the
action. It is straightforward to verify that $d^D_n =
(\partial_{n-1})_*$ induces a map on coinvariants, yielding a
complex $\ov{D}_* = \{\ov{D}_n, d^{\ov{D}}_n\}$. Let
$\{E_n\}_{n\ge 0}$ be the chain complex $E_n =
H_1(B\Gamma_{n-1};\Bbb Q)$, $d^E_n =
\sum_{i=0}^{n-1}(-1)^i(\partial_i)_* : E_n\to E_{n-1}$, with $E_n
= 0$ for $n < 1$. Finally, let $\ov{E}_n$ denote the quotient of
$E_n$ by the subgroup generated by elements of the form
$(s_j)_*(x)$ where $x\in H_1(B\Gamma_{n-2};\Bbb Q)$ and $(s_j)_*$
the map induced on homology by the degeneracy $s_j$. Again,
$d^E_n$ descends to a map $d^{\ov{E}}_n:\ov{E}_n\to \ov{E}_{n-1}$,
making $\ov{E}_* = \{\ov{E}_n, d^{\ov{E}}_n\}$ a chain complex.

\proclaim{\bf\underbar{Lemma 5.10}} There are isomorphisms
$$
H_*(\ov{D}_*)\cong H_*(E_*)\cong H_*(\ov{E}_*)\cong H_*(B\pi;\Bbb Q)
$$
\endproclaim

\prf First, $\ov{E}_*$ is the normalized complex of $E_*$ formed
by collapsing the acyclic subcomplex of degenerate elements, so
$H_*(E_*)\cong H_*(\ov{E}_*)$. The $E^1$-term of the spectral
sequence in homology $\{E^1_{p,q}= H_p(B\Gamma_q;\Bbb Q)\}$
converging to  $H_*(B\Gamma\hskip-.03in .;\Bbb Q)\cong
H_*(B\pi;\Bbb Q)$ vanishes for $p > 1$ (as $\Gamma_n$ is free for
all $n\ge 0$); upon passing to the $E^2$-term, the $p=0$ line
vanishes while the $p=1$ line yields precisely the homology of
$E_*$. Thus $H_*(E_*)\cong H_*(B\pi;\Bbb Q)$.

To see the relationship between $H_*(\ov{D})$ and $H_*(\ov{E}_*)$,
we construct maps in both directions. Let $\lambda^j_n =
s_j\partial_j : \Gamma_n\to \Gamma_n$ where $0\le j < n$.
Inductively define $r^j_n(g)$ by $r^{-1}_n(g) = g;\, r^j_n(g) =
(r^{j-1}_n(g))(\lambda^j_n(r^{j-1}_n(g)))^{-1}$. For each $j < n$,
$r^j_n$ defines a projection of sets $r^j_n:\Gamma_n\surj
\Gamma^j_n$. Let $\Gamma^{j,d}_n$ denote the normal subgroup of
$\Gamma^j_n$ generated by commutators of the form $[x,y]$ where
$x\in \Gamma^j_n$ and $y=s_i(z)$ for some $z\in \Gamma_{n-1}$ and
$0\le i\le j$. Let $\ov{\Gamma}^j_n = \Gamma^j_n/\Gamma^{j,d}_n$.
An application of simplicial identities verifies that the
composition $\Gamma_n\overset r^j_n\to\longrightarrow
\Gamma^j_n\surj \ov{\Gamma}^j_n$ is a homomorphism. There is an
equality $H_1(B\Gamma^{n-1}_n;\Bbb Q)_{\Gamma_n} =
H_1(B\ov{\Gamma}^{n-1}_n;\Bbb Q)$, from which we see that
$r^{n-2}_{n-1}$ induces for each $n\ge 1$ a homomorphism
$(r^{n-2}_{n-1})_*:E_n\to\ov{D}_n$. Again, use of simplicial
identities allows one to work out the following recursive
description for $\partial_{n}(r^{n-1}_{n}(g))$: let $A_0(g) =
r^{n-2}_{n-1}(\partial_0(g))$, with  $A_i(g) =
r^{n-2}_{n-1}(\partial_i(g))A_{i-1}(g)^{-1}$. Then
$\partial_n(r^{n-1}_n(g)) = A_n(g)$. This formula implies
$\{(r^{n-2}_{n-1})_*\}_{n\ge 1}$ defines a chain map $E_*\to
\ov{D}_*$. On the other hand, the inclusion
$\Gamma^{n-1}_n\hookrightarrow \Gamma_n$ for each $n$ defines a
chain map $D_*\to E_*$ which factors by the projection $D_*\surj
\ov{D}_*$, as $H_1(B\Gamma_n;\Bbb Q) = H_1(B\Gamma_n;\Bbb
Q)_{\Gamma_n}$. It is then straightforward to check $\ov{D}_*\to
E_*\to \ov{D}_*$ is the identity, and that the composition $E_*\to
\ov{D}_*\to E_*\surj \ov{E}_*$ agrees with the projection
$E_*\surj \ov{E}_*$. Together these yield the isomorphism
$H_*(\ov{D})\cong H_*(E_*)$.\hfill //
\endpf
\vskip.2in

\underbar{Remark 5.11} The surjection $D_*\surj \ov{D}_*$ in
general yields neither a surjection nor injection upon passage to
homology; in fact it is not hard to work out an explicit
description of the homology groups themselves which indicates the
difference. From the exactness of the complex
$\{\Gamma^{n-2}_{n-1},\partial_{n-1}\}_{n\ge 1}$ we see that
$H_1(D_*) = H_1(\ov{D}_*)$, while for $n \ge 2$ one has
$$
\gathered
H_n(D_*) =
\frac{\Gamma^{n-2}_{n-2}\cap [\Gamma^{n-3}_{n-2},\Gamma^{n-3}_{n-2}]}{[\Gamma^{n-2}_{n-2},\Gamma^{n-2}_{n-2}]}\\
H_n(\ov{D}_*) =
\frac{\Gamma^{n-2}_{n-2}\cap [\Gamma^{n-3}_{n-2},\Gamma_{n-2}]}{[\Gamma^{n-2}_{n-2},\Gamma_{n-2}]}
\endgathered
$$
These groups are the same for $n = 2$ but will differ
in general when $n>2$. When $n=2$, the expression for $H_2(D_*) =
H_2(\ov{D}_*)$ is precisely Hopf's formula for $H_2(B\pi)$, so the
description of $H_n(\ov{D}_*)$ when $n>2$ may be seen as a higher
dimensional analogue of this formula. The equalities can be
formulated integrally, and hold integrally. \vskip.2in

By Lemma 5.10, we see that any $0\ne x\in H_n(B\pi;\Bbb Q) =
H_n(\ov{D}_*)$ maps non-trivially to
$\ov{D}_n/im(d_{n+1}^{\ov{D}})$, and in turn can be lifted
(non-uniquely) to an element $0\ne y_x\in D_n/im(d_{n+1}^{D})$
under the surjection $D_n/im(d_{n+1}^{D})\surj
\ov{D}_n/im(d_{n+1}^{\ov{D}})$. Consider the following commutative
diagrams (in which $A\otimes\Bbb Q$ is abbreviated by $A^{\Bbb
Q}$).
$$
\diagram
& &H_n(B\pi;\Bbb Q)\dto|<<\tip\\
H_1(B\Gamma^{n-1}_n;\Bbb Q)_{\Gamma_n}\rto^{(\ov{\partial}_n)_*}
& H_1(B\Gamma^{n-2}_{n-1};\Bbb Q)_{\Gamma_{n-1}}\rto
& coker_1^{\ov{D},n}\\
H_1(B\Gamma^{n-1}_n;\Bbb Q)\rto^{(\partial_n)_*}\dto\uto|>>\tip
& H_1(B\Gamma^{n-2}_{n-1};\Bbb Q)\rto\dto\uto|>>\tip
& coker_1^{D,n}\dto\uto|>>\tip\\
K_1^a(\partial_n(I[\Gamma^{n-1}_n]))^{\Bbb Q}\rto\dto
&K_1^a(I[\Gamma^{n-2}_{n-1}])^{\Bbb Q}\rto\dto
& coker_2^n\dto\\
K_1^a(\partial_n(C^i(\Gamma\hskip-.03in .^+)^{n-1}_n))^{\Bbb Q}\rto\dto
&K_1^a(C^i(\Gamma\hskip-.03in .^+)^{n-2}_{n-1})^{\Bbb Q}\rto\dto
& coker_3^{a,n}\dto\\
K_1^t(\partial_n(C^i(\Gamma\hskip-.03in .^+)^{n-1}_n))^{\Bbb
Q}\rto &K_1^t(C^i(\Gamma\hskip-.03in .^+)^{n-2}_{n-1})^{\Bbb
Q}\rto & coker_3^{t,n}
\enddiagram
\tag5.12
$$
Here $I[G] = ker(\Bbb Z[G]\to\Bbb Z)$ and
$H_1(B\Gamma^{n-2}_{n-1};\Bbb Q)\to
K_1^a(I[\Gamma^{n-2}_{n-1}])^{\Bbb Q}$ is the (reduced,
restricted) assembly map for $\Gamma^{n-2}_{n-1}$ in algebraic
$K$-theory. The map $H_1(B\Gamma^{n-1}_{n};\Bbb Q)\to
K_1^a(\partial_n(I[\Gamma^{n-1}_{n}]))^{\Bbb Q}$ is the
composition $H_1(B\Gamma^{n-1}_{n};\Bbb Q)\to
K_1^a(I[\Gamma^{n-1}_{n}])^{\Bbb Q}\to
K_1^a(\partial_n(I[\Gamma^{n-1}_{n}]))^{\Bbb Q}$. The vertical
maps from the fourth to the fifth lines are induced by the natural
transformation $\Bbb Z[\,\,\,]\to C^i(\,\,\,)$; from the fifth to
the sixth lines by the transformation from algebraic to
topological $K$-theory. Finally, the groups on the right are the
cokernels of the respective maps on each line, while the vertical
maps are those induced by maps of cokernels.

Let $[c]\in H^n(B\pi;\Bbb Q)$. As noted above in (5.5), $[c]$
identifies with a simplicial homomorphism
$\phi_c:\Gamma\hskip-.03in .^+\to K(\Bbb Q,n-1)\hskip-.015in .$;
in turn $\phi_c$ induces a simplicial algebra homomorphism
$C^i(\Gamma\hskip-.03in .^+)\to C^i(K(\Bbb Q,n-1)\hskip-.015in .)$
which we also denote be $\phi_c$. Following (5.9) our second
diagram is
$$
\diagram H_1(B\Gamma^{n-1}_n;\Bbb Q)\rto^{(\partial_n)_*}\dto &
H_1(B\Gamma^{n-2}_{n-1};\Bbb Q)\rto\dto
& coker_1^{D,n}\dto\\
K_1^a(\partial_n(C^i(\Gamma\hskip-.03in .^+)^{n-1}_n))^{\Bbb
Q}\rto\dto &K_1^a(C^i(\Gamma\hskip-.03in .^+)^{n-2}_{n-1})^{\Bbb
Q}\rto\dto
& coker_3^{a,n}\dto\\
K_1^a(\partial_n(C^i(K(\Bbb Q,n-1)\hskip-.015in .)^{n-1}_n))^{\Bbb
Q}\rto\ddouble &K_1^a(C^i(K(\Bbb Q,n-1)\hskip-.015in
.)^{n-2}_{n-1})^{\Bbb Q}\rto\ddouble
& coker_4^{a,n}\dto\\
K_1^a(I^i(\Bbb Q)^2)^{\Bbb Q}\rto & K_1^a(I^i(\Bbb Q))^{\Bbb
Q}\rto & K_1^a(I^i(\Bbb Q)/I^i(\Bbb Q)^2)^{\Bbb Q}\dto^{\cong}\\
& & \Bbb C
\enddiagram
\tag5.13
$$

Given $0\ne x\in H_n(B\pi;\Bbb Q)$, we choose a lift $y_x\in
coker_1^{D,n}$ which maps to the image of $x$ in
$coker_1^{\ov{D},n}$. Mapping $y_x$ by the composition on the
right-hand side of (5.13) yields an element in $\Bbb C$ equal to
the image of $<[c],x>\in\Bbb Q$ under the inclusion $\Bbb
Q\hookrightarrow \Bbb C$. In particular, it depends only on $x$
and not the choice of lift $y_x$. Referring to these diagrams, we
may summarize this as

\proclaim{\bf\underbar{Theorem 5.14}} The Strong Novikov
Conjecture is equivalent to the statement:

(i) For all $n\ge 1$ and $y\in coker_1^{D,n}$, if $y$ projects to
a non-zero element in $im(H_n(B\pi;\Bbb Q)\inj
coker_1^{\ov{D},n})$, then it maps to a non-zero element in
$coker_3^{t,n}$.

Moreover, the algebraic analogue of (i) is always true:

(ii) For all $n\ge 1$ and $y\in coker_1^{D,n}$, if $y$ projects to
a non-zero element in $im(H_n(B\pi;\Bbb Q)\inj
coker_1^{\ov{D},n})$, then it maps to a non-zero element in
$coker_3^{a,n}$.
\endproclaim

\prf The first part is a restatement of Theorem 5.4, in the
context of (5.12). The second part follows from the Universal
Coefficient Theorem.\hfill //
\endpf

This method of detecting the image of the restricted assembly map does not work for
topological $K$-theory directly, as the inclusion
$(I^i(\Bbb Q))^2\hookrightarrow I^i(\Bbb Q)$ induces an isomorphism of topological
$K$-groups.

Theorem 5.14 can be restated in terms of absolute $K$-groups. Let
$\wt C^i_{n-1}(\Gamma\hskip-.03in .^+) = C^i(\Gamma\hskip-.03in
.^+)^{n-2}_{n-1}/\partial_n(C^i(\Gamma\hskip-.03in
.^+)^{n-1}_{n})$ equipped with the fine topology.

\proclaim{\bf\underbar{Theorem 5.15}} The Strong Novikov
Conjecture is equivalent to the statement:

(i) For all $n\ge 1$ and $y\in coker_1^{D,n}$, if $y$ projects to
a non-zero element in $im(H_n(B\pi;\Bbb Q)\inj
coker_1^{\ov{D},n})$, then it maps to a non-zero element in
$K_1^t(\wt C^i_{n-1}(\Gamma\hskip-.03in .^+))^{\Bbb Q}$.

Moreover, the algebraic analogue of (i) is always true:

(ii) For all $n\ge 1$ and $y\in coker_1^{D,n}$, if $y$ projects to
a non-zero element in $im(H_n(B\pi;\Bbb Q)\inj
coker_1^{\ov{D},n})$, then it maps to a non-zero element in
$K_1^a(\wt C^i_{n-1}(\Gamma\hskip-.03in .^+))^{\Bbb Q}$.
\endproclaim

These methods also apply to the rational Witt groups of the group
algebra. The augmented simplicial group algebra $\Bbb
Z[\Gamma\hskip-.03in .^+]$ is a resolution of $\Bbb Z[\pi]$ if
$\Gamma\hskip-.03in .^+$ is a resolution of $\pi$, yielding a
short-exact sequence of rings-with-involution
$$
\Bbb Z[\Gamma\hskip-.03in .^+]^n_n
\inj \Bbb Z[\Gamma\hskip-.03in .^+]^{n-1}_n
\surj \Bbb Z[\Gamma\hskip-.03in .^+]^{n-1}_{n-1}\qquad\qquad n\ge 0
$$
As $W_*(_-)^{\Bbb Q}$ satisfies excision, there are boundary maps
$\partial_*: W_m(\Bbb Z[\Gamma\hskip-.03in .^+]^{n-1}_{n-1})^{\Bbb
Q}\to W_{m-1}(\Bbb Z[\Gamma\hskip-.03in .^+]^{n}_{n})^{\Bbb Q}$.
Composing them yields a map of Witt groups $$
\partial^{(n-1)}: W_n(\Bbb Z[\pi])^{\Bbb Q}\to W_1(\Bbb Z[\pi]^{n-2}_{n-2})^{\Bbb Q}
$$
for each $n\ge 1$. Using the periodicity of the rationalized Witt
groups and fact that the assembly map $H_*(BG;\Bbb Q)\otimes
W_*(\Bbb Z)^{\Bbb Q}\to W_*(\Bbb Z[G])^{\Bbb Q}$ is an isomorphism
when $G$ is a free group, the same line of argument as above
implies

\proclaim{\bf\underbar{Theorem 5.16}} The Novikov conjecture is
true for $\pi$ if and only if the composition $$ H_n(B\pi;\Bbb
Q)\to W_n(\Bbb Z[\pi])^{\Bbb
Q}\overset\partial^{(n-1)}\to\longrightarrow W_1(\Bbb
Z[\pi]^{n-2}_{n-2})^{\Bbb Q}
$$
is injective for each $n\ge 1$.
\endproclaim

The first map in this composition is induced by the inclusion of
the summand $H_*(B\pi;\Bbb Q)\inj H_*(BG;\Bbb Q)\otimes W_*(\Bbb
Z)^{\Bbb Q}$, followed by the assembly map to $W_*(\Bbb
Z[G])^{\Bbb Q}$.
\vskip.3in


\head 6 Some additional remarks
\endhead
\vskip.2in

We conclude by describing some additional properties of $C^i$-algebras.
As we have noted, $C^i(_-)$ is functorial with
respect to group homomorphisms, hence extends to a functor from $(s.gps)$ (the category of
simplicial groups) to $(s.t.alg.)$ (the category of simplicial topological algebras). Recall
also that a simplicial set $X\hskip-.015in .$ is said to be $r$-reduced if
$X_k = \{pt\}$ for all $k\le r$.

\proclaim{\bf\underbar{Lemma 6.1}} Let $G\hskip-.015in .$ be a
$0$-reduced simplicial group. Then the inclusion of simplicial
algebras $\Bbb C[G\hskip-.015in .]\hookrightarrow
C^i(G\hskip-.015in .)$ is a weak homotopy equivalence (the
homotopy groups being the simplicial homotopy groups as defined in
section 1).
\endproclaim

\prf Consider the commuting diagram
$$
\diagram
\Bbb C[G\hskip-.015in .]\rto\dto & \Bbb C[G\hskip-.015in .]\,\Hat{}\dto\\
C^i(G\hskip-.015in .)\rto & C^i(G\hskip-.015in .)\,\Hat{}
\enddiagram
$$
where the horizontal arrows denote completion with respect to the (simplicial) augmentation
ideal, while the vertical arrows arise from the natural transformation
$\Bbb C[_-]\hookrightarrow C^i(_-)$. By Curtis convergence [Cu], the horizontal maps are weak equivalences. By
Lemma 5.7, the right vertical map is an isomorphism. The result follows.\hfill //
\endpf

For a simplicial algebra $A\hskip-.015in .$ over $\Bbb C$,
$C_*(A\hskip-.015in .)$ resp. $CC_*(A\hskip-.015in .)$ will denote
the total Hochschild resp. cyclic complex associated to the
simplicial complexes $\{[n]\mapsto C_*(A_n)\}_{n\ge 0}$ resp.
$\{[n]\mapsto CC_*(A_n)\}_{n\ge 0}$. For algebras topologized by
the fine topology, the algebraic and topological complexes are the
same, so there is no need to distinguish between the two.
\vskip.1in

\underbar{Question 6.2} Is $C^i(\Gamma\hskip-.03in .)$ a
resolution of $\wt C^i(\Gamma\hskip-.03in .) \overset def\to =
\pi_0(C^i(\Gamma\hskip-.03in .))$? \vskip.1in

\underbar{Question 6.3} Does the image of $H_*(B\pi;\Bbb Q)\to
HH_*(\Bbb C[\pi])\to HH_*(\wt C^i(\Gamma\hskip-.03in .))$ map
injectively under $I:HH_*(\wt C^i(\Gamma\hskip-.03in .))\to
HC_*(\wt C^i(\Gamma\hskip-.03in .))$? \vskip.1in

\proclaim{\bf\underbar{Theorem 6.4}} If the answer to these two
questions is \lq\lq yes\rq\rq for a given resolution
$\Gamma\hskip-.03in .$ of $\pi$, then the Strong Novikov
Conjecture is true for $\pi$.
\endproclaim

\prf If $C^i(\Gamma\hskip-.03in .)$ is a resolution of $\wt
C^i(\Gamma\hskip-.03in .)$, then the augmentation
$C^i(\Gamma\hskip-.03in .)\surj \wt C^i(\Gamma\hskip-.03in .)$
induces a quasi-isomorphism of chain complexes
$C_*(C^i(\Gamma\hskip-.03in .))\overset\simeq\to\surj C_*(\wt
C^i(\Gamma\hskip-.03in .))$. Given $0\ne x\in H_n(B\pi;\Bbb Q)$,
choose a cohomology class $[\phi]\in H^n(B\pi;\Bbb Q)$ which pairs
non-trivially with $x$, and then a simplicial group homomorphism
$\phi :\Gamma\hskip-.03in .\to K(\Bbb Q,n-1)\hskip-.015in .$
representing $[\phi]$ in the manner discussed above. We may assume
without loss of generality that $n
> 1$, as injectivity in the case $n = 1$ may be verified directly.
Now consideration of the commuting diagram
$$
\diagram
H_*(B\Gamma\hskip-.03in .;\Bbb Q)\rto\dto^{\cong} & HH_*(C^i(\Gamma\hskip-.03in .))
\rto^(.4){(\phi)_*}\dto^{\cong} & HH_*(C^i(K(\Bbb Q,n-1)\hskip-.015in .))\dto^{\cong}\\
H_*(B\pi;\Bbb Q)\rto & HH_*(\wt C^i(\Gamma\hskip-.03in .)) & HH_*(\Bbb C[K(\Bbb Q,n-1)\hskip-.015in .])
\enddiagram
$$
verifies injectivity of the composition
$H_*(B\pi;\Bbb Q)\to HH_*(\Bbb C[\pi])\to HH_*(\wt C^i(\Gamma\hskip-.03in .))$. If the
composition of this map with $I: HH_*(\wt C^i(\Gamma\hskip-.03in .))\to
HC_*(\wt C^i(\Gamma\hskip-.03in .))$ is still injective, as stated in (6.3), then the
restricted assembly map
$$
H_*(B\pi;\Bbb Q)\to K_*^t(\wt C^i(\Gamma\hskip-.03in .))
$$
is injective, via the Chern-Connes-Karoubi character defined in [O1, App.]. Finally,
by a five lemma argument one verifies that the natural map
$\wt C^i(\Gamma\hskip-.03in .)\to C^*(\pi)$ induces an isomorphism of topological $K$-groups
in dimensions $\ge 1$. By what we have previously shown, this implies the Strong Novikov
Conjecture for $\pi$.\hfill //
\endpf

It is not hard to show that given two free simplicial
resolutions $\Gamma\hskip-.03in .$, $\Gamma\hskip-.03in .'$ of
$\pi$, the resulting simplicial algebras $C^i(\Gamma\hskip-.03in
.)$ and $C^i(\Gamma\hskip-.03in .')$ are homotopy equivalent. Thus
the homotopy groups of $C^i(\Gamma\hskip-.03in .)$ are invariants
of $\pi$, and either the answer to (6.2) is \lq\lq yes\rq\rq or
one has a new set of invariants associated to the discrete group
$\pi$.


\enddocument